\documentclass{amsart}
\newtheorem{prop}{Proposition}
\newtheorem{lemma}{Lemma}

\newtheorem{theo}{Theorem}
\theoremstyle{definition}
\newtheorem{exam}{Example}
\newtheorem{rem}{Remark}

\newcommand{\myA}{\mathcal{A}_q}
\newcommand{\myC}{\mathbb{C}}
\newcommand{\myI}{I}
\newcommand{\myL}{\mathrm{Li}_2}
\newcommand{\myNP}{\mathbb{Z}_{\ge0}}
\newcommand{\myP}{\phi}

\newcommand{\myRp}{\mathbb{R}_{\ge 0}}
\newcommand{\myS}{S}
\newcommand{\myU}{\mathsf{u}}
\newcommand{\myV}{\mathsf{v}}
\newcommand{\myZ}{\mathbb{Z}}
\newcommand{\myZp}{\mathbb{Z}_{\ge0}}

\title{The $q$-binomial formula and the Rogers dilogarithm identity}
\author{R.M.  Kashaev}
\address{Universit\'e de Gen\`eve,
Section de math\'ematiques,
2-4, rue du Li\`evre,
CP 240,
1211 Gen\`eve 24, Suisse
}
\address{V.A. Steklov Institute of Mathematics at St. Petersburg , 27, 
Fontanka, St. Petersburg 191023, Russia}
\email{Rinat.Kashaev@math.unige.ch}
\date{July 2004}
\thanks{The work is supported in part by the Swiss National Science Foundation}
\begin{document}

\begin{abstract}
The $q$-binomial formula in the limit $q\to 1^-$ is shown to be equivalent to the Rogers five term dilogarithm identity.
\end{abstract}
\maketitle
\section{Introduction}

For any $q,x\in]0,1[$ define a $q$-exponential function as an infinite 
product
\[
\myP(x):=1/(x;q)_\infty :=1/\prod_{n\ge0}(1-q^n x)
\]
A finite product
\[
(x;q)_k:=\prod_{n=0}^{k-1}(1-q^n x),\quad \forall k\in\myZp
\]
can be expressed as a ratio of two $q$-exponentials:
\[
(x;q)_k=\frac{(x;q)_\infty}{(xq^k;q)_\infty}=\frac{\myP(xq^k)}{\myP(x)}
\]
The $q$-binomial formula (see, for example, \cite{GR}) 
is given by the following identity
\begin{equation}\label{e15}
\sum_{n\ge 0}\frac{(a;q)_n}{(q;q)_n}z^n=\frac{(az;q)_\infty}{(z;q)_\infty},\quad |z|<1
\end{equation}
which, by using the above notation, can also be written entirely 
in terms of the function $\myP(x)$:
\begin{equation}\label{e1}
\sum_{n\ge 0}
\frac{\myP(aq^{n})}{\myP(q^{n+1})}z^n=\frac{\myP(a)\myP(z)}{\myP(q)\myP(az)}
\end{equation}
 The following expansion formulas
\begin{equation}\label{e10}
\myP(x)=\sum_{n\ge0}\frac{x^n}{(q;q)_n}
\end{equation}
and
\begin{equation}\label{e11}
\frac1{\myP(x)}=\sum_{n\ge0}\frac{(-1)^n q^{n(n-1)/2}x^n}{(q;q)_n}
\end{equation}
are both particular cases of the $q$-binomial formula.

The asymptotic formula 
\[
\myP(x)\sim e^{-\myL(x)/\ln q},\quad q\to 1^-
\] 
where
\[
\myL(x):=\sum_{n=1}^\infty \frac{x^n}{n^2}
\]
is the Euler dilogarithm function, has been used in \cite{FK,Kir}
to give an interpretation to $\myP(x)$ as a quantum version of the
dilogarithm function. 
In particular, by using a formal reasoning coming from quantum mechanics,
 it has been shown that the quantum five term identity
\begin{equation}\label{e13}
\myP(\myU)\myP(\myV)=\myP(\myV)\myP(-\myV\myU)\myP(\myU)
\end{equation}
where $\myP(\myU)$, $\myP(\myV)$, and $\myP(-\myV\myU)$ are elements
in the algebra $\myA=\myC_q[[\myU,\myV]]$ of formal power series in two 
elements $\myU, \myV$ satisfying the commutation relation 
$\myU\myV=q\myV\myU$, in the limit $q\to 1^-$ reproduces the Rogers
pentagonal identity for the dilogarithm
\begin{equation}\label{e16}
\myL(a)+\myL(z)=\myL(az)+\myL\left(\frac{a-az}{1-az}\right)
+\myL\left(\frac{z-az}{1-az}\right)
+\ln\left(\frac{1-z}{1-az}\right)\ln\left(\frac{1-a}{1-az}\right)
\end{equation}
The purpose of this paper is to make 
the statement of the paper \cite{FK} mathematically rigorous\footnote{I would like to thank Yu. Manin for posing this question.}. Namely, we first show that the
 identity~\eqref{e13} is related to the $q$-binomial
formula~\eqref{e15} and then derive from the latter the Rogers
identity~\eqref{e16} in the limit $q\to 1^-$. 
 The main result follows.
\begin{theo}\label{t1}
Let $q,a,z\in]0,1[$. Then in the limit $q\to 1^-$ the $q$-binomial
identity~\eqref{e1} 
leads to the Rogers pentagonal identity~\eqref{e16}.
\end{theo}
The rest of this paper is organized as follows. In Section~\ref{sec2} the equivalence between he $q$-binomial formula and the quantum pentagonal identity is explained, while Section~\ref{sec3} contains the proof of Theorem~\ref{t1}.
 
\section{The $q$-binomial formula and the quantum pentagonal identity}\label{sec2}

The relation between the formulas~\eqref{e15} and \eqref{e13} can be
established by comparing the expansion coefficients of $a^m z^n$ in
 \eqref{e15} and $\myV^n\myU^m$ in \eqref{e13}, respectively. 
\begin{prop}
The $q$-binomial formula is equivalent to the following set of identities
\begin{equation}\label{e12}
\frac{q^{mn}}{(q;q)_m(q;q)_n}
=\sum_{k=0}^{\min(m,n)}\frac{(-1)^k q^{k(k-1)/2}}{(q;q)_{m-k}(q;q)_{n-k}(q;q)_k},\quad \forall m,n\in\myNP
\end{equation}
\end{prop}
\begin{proof}
Let us write the  $q$-binomial formula in the form
\[
\sum_{n\ge 0}\frac{\myP(aq^n)}{(q;q)_n}z^n=\frac{\myP(a)\myP(z)}{\myP(az)}
\]
or, using formula~\eqref{e10} in the left hand side, we have
\[
\sum_{m,n\ge 0}\frac {q^{mn}a^mz^n}{(q;q)_m(q;q)_n}=\frac{\myP(a)\myP(z)}{\myP(az)}
\]
Again, using the expansion formulas~\eqref{e10}, \eqref{e11} in the
right hand side,  and equating the coefficients of the monomials
$a^mz^n$ in both sides of the equality, we arrive at formula~\eqref{e12}.
\end{proof}
\begin{prop}
The set of identities~\eqref{e12} is equivalent to the
quantum five term identity \eqref{e13}. 
\end{prop}
\begin{proof}
We multiply the both sides of \eqref{e12} by $\myV^n\myU^m$ and sum
over $m$ and $n$. The result 
can be easily written in the form of equation~\eqref{e13} by using the commutation relation $\myU\myV=q\myV\myU$, and, in particular, the formula $\myV^k\myU^kq^{k(k-1)/2}=(\myV\myU)^k$.
\end{proof}

\section{Proof of Theorem~\ref{t1}}\label{sec3}

\begin{lemma}\label{l2}
Let $k,l\in\myZ$ be such that $k\le l$ and $f_\pm\colon [k,l+1]\to\myRp $ be functions, where $f_-$ is decreasing and $f_+$ is increasing. Then 
\begin{equation}\label{e17}
\sum_{n=k+1}^{l+1}f_-(n)\le\int_{k}^{l+1} f_-(t)dt\le\sum_{n=k}^lf_-(n)
\end{equation}
\begin{equation}\label{e171}
\sum_{n=k}^{l}f_+(n)\le\int_{k}^{l+1} f_+(t)dt\le\sum_{n=k+1}^{l+1}f_+(n)
\end{equation}
\end{lemma}
\begin{proof}
The inequality 
\[
f_-(n+1)\le f_-(x)\le f_-(n),\quad \forall n\in\myZ\cap [k,l],\ \forall x\in[n,n+1]
\]
implies that $f_-(n+1)\le \int_n^{n+1}f_-(x)dx\le f_-(n)$. Thus, summing over all possible $n$ we arrive
at formula~\eqref{e17}. The proof of formula~\eqref{e171} is similar.
\end{proof}
\begin{rem}
The variables $k$ and $l$ in Lemma~\ref{l2} can take infinite values $k=-\infty$ or $l=\infty$.
\end{rem}
In what follows, for any function $f\colon \myRp\to\myRp$ we shall use the notation
\[
\myS(f):=\sum_{n\ge 0}f(n),\quad \myI(f):=\int_{0}^\infty f(t)dt
\]
If a decreasing function $f\colon \myRp\to\myRp$ is integrable on $\myRp$ then, as a particular case of Lemma~\ref{l2}, we have
\[
\myS(f)-f(0)\le\myI(f)\le\myS(f)
\] 
or equivalently
\begin{equation}\label{e18}
0\le \myS(f)-\myI(f)\le f(0)
\end{equation}
\begin{exam}
The function $f(t)=-\ln(1-q^tx)$ is decreasing and integrable on $\myRp$, and 
\[
\myS(f)=\ln\myP(x),\quad \myI(f)=-\int_{0}^\infty \ln(1-q^tx)dt=
\frac1 {\ln q}\int_0^x\ln(1-z)\frac{dz}{z}=-\frac{\myL(x)}{\ln q}
\]
Thus, for any $q,x\in]0,1[$ inequalities~\eqref{e18} imply that 
\begin{equation}\label{e2}
1\le \myP(x) e^{\myL(x)/\ln q}\le \frac1{1-x}
\end{equation}
\end{exam}
\begin{lemma}\label{l1}
Let $g\colon \myRp\to\myRp $ be an integrable function increasing
in the segment $[0,x_0]$ and  decreasing on the interval $[x_0,\infty[$. Let
also $n_0\in\myNP$ be such that $g(n)\le g(n_0)$ for all $n\in\myNP$ 
($n_0$ is equal either to $[x_0]$ (the integer part of $x_0$) or
$[x_0]+1$). Then
\begin{equation}\label{e172}
g(n_0)\le\sum_{n\ge0}g(n)\le \int_0^\infty g(x)dx+g(n_0)
\end{equation}
\end{lemma}
\begin{proof}
The inequality $g(n_0)\le\sum_{n\ge0}g(n)$ follows directly from the positivity of $g(x)$. To prove the second part of \eqref{e172}, note that we can apply Lemma~\ref{l2} to functions $f_+=g\vert_{[0,[x_0]]}$
and $f_-=g\vert_{[[x_0]+1,\infty[}$. Thus, the left hand sides of the inequalities in Lemma~\ref{l2} take the forms
\[
\sum_{n=0}^{[x_0]-1}g(n)\le \int_0^{[x_0]}g(x)dx,\quad
\sum_{n=[x_0]+2}^\infty g(n)\le \int_{[x_0]+1}^{\infty}g(x)dx
\]
Adding these to each other, we obtain
\[
\sum_{n=0}^{\infty}g(n)-g([x_0])-g([x_0]+1)\le\int_0^{\infty}g(x)dx-\int_{[x_0]}^{[x_0]+1}g(x)dx
\]
which, combined with the inequality
\[
\int_{[x_0]}^{[x_0]+1}g(x)dx\ge g(n_0')
\]
where $\{n_0,n_0'\}=\{[x_0],[x_0]+1\}$, is equivalent to the second part of \eqref{e172}.
\end{proof}
\begin{prop}
There exists $\epsilon\in]0,1[$ such that for any 
$q\in]1-\epsilon,1[$
the function
\[
g(x)=\frac{\myP(a q^{x})} {\myP(q^{1+x})}z^x
\]
where $a,z\in]0,1[$,
 satisfies the conditions of Lemma~\ref{l1}. 
 \end{prop}
 \begin{proof}
 The integrability of $g(x)$ is evident.
We have the following formula for its derivative
\[
\frac{g'(x)}{g(x)}=\ln z 
-\ln (q)(q-a)\myS(h_x)
\]
where
\[
h_x(t)=\frac{q^{x+t}}{(1-q^{1+x+t})
(1-a q^{x+t})}
\]
satisfies the conditions of Lemma~\ref{l2} so that
\[
\myS(h_x)\ge\myI(h_x)=-\frac1{\ln(q)(q-a)}\ln\left(\frac{z(1-aq^x)}{1-q^{1+x}}\right)
\]
Evidently, the function $\myS(h_x)$ is decreasing in $x$.
Assuming that $q>1-z(1-a)$, we obtain
\[
\frac{g'(0)}{g(0)}\ge \ln\left(\frac{z(1-a)}{1-q}\right)>0
\]
Besides, it is easy to see that
\[
\lim_{x\to\infty}\frac{g'(x)}{g(x)}=\ln z <0
\]
Thus, we have shown that for $\epsilon=z(1-a)$ and any $q\in]1-\epsilon,1[$ the continuous function $g'(x)/g(x)$
is decreasing, positive at $x=0$ and negative for sufficiently large $x$, i.e. there exists unique
$x_0\in]0,\infty[$ such that $g'(x_0)=0$ and all conditions of Lemma~\ref{l1} are satisfied.
\end{proof}
\begin{prop}\label{p3}
\begin{equation}\label{e3}
\lim_{q\to 1^-}\ln(q)\ln\myS(g)=F(\xi_0),\quad\xi_0=\frac{1-z}{1-az}
\end{equation}
where
\[
F(\xi)=\myL(\xi)-\myL(a\xi)+\ln(\xi)\ln(z)
\] 
\end{prop}
\begin{proof}
For any $\xi\in]0,1[$ equation~\eqref{e2} implies that
\[
\lim_{q\to 1^-}\ln(q)\ln( g\left({\ln\xi}/{\ln q}\right))=F(\xi)
\]
Thus, one has asymptotically
\[
g(\left({\ln\xi}/{\ln q}\right)\sim e^{\frac{F(\xi)}{\ln q}},\quad q\to 1^-
\]
and, by using the steepest decent method, one has also
\[
\myI(g)\sim e^{\frac{F(\xi_0)}{\ln q}},\quad q\to 1^-
\]
where $\xi_0={(1-z)}/{(1-az)}\in]0,1[$ is the unique solution of the equation $F'(\xi)=0$.
The asymptotic formula for $S(g)$ follows immediately from Lemma~\ref{l1} after taking into account the fact that 
\[
x_0\sim n_0\sim \frac{\ln\xi_0}{\ln q},\quad q\to 1^-
\]
and, correspondingly, $g(n_0)\sim g(x_0)\sim\myI(g)$, $q\to 1^-$.
\end{proof}
\begin{proof}[Proof of Theorem~\ref{t1}]
Using Lemma~\ref{l2}, we have immediately
\[
\lim_{q\to 1^-}\ln(q)\ln\left(\frac{\myP(a)\myP(z)}{\myP(q)\myP(az)}\right)=
\myL(1)+\myL(az)-\myL(a)-\myL(y)
\]
Combining this formula with equation~\eqref{e3}, we conclude that the $q$-binomial
identity~\eqref{e1} leads to the following identity:
\[
F(\xi_0)=\myL(1)+\myL(az)-\myL(a)-\myL(z)
\]
or explicitly,
\[
\myL(\xi_0)-\myL(a\xi_0)+\ln(\xi_0)\ln(z)=\myL(1)+\myL(az)-\myL(a)-\myL(z)
\]
which we rewrite in the form
\[
\myL(a)+\myL(z)=\myL(az)+\myL(a\xi_0)+\myL(1)-\myL(\xi_0)-\ln(\xi_0)\ln(z)
\]
Using the identity
\[
\myL(x)+\myL(1-x)=\myL(1)-\ln(x)\ln(1-x),\quad \forall x\in[0,1]
\]
we rewrite it further
\[
\myL(a)+\myL(z)=\myL(az)+\myL(a\xi_0)+\myL(1-\xi_0)+\ln(\xi_0)\ln\left((1-\xi_0)/z\right)
\]
which is exactly  the Rogers identity~\eqref{e16}.

\end{proof}

\end{document}